\documentclass[12pt]{amsart}
\usepackage{amsmath}
\usepackage{amsthm}
\usepackage{amscd}

\newcommand{\Q}{{\mathbb Q}}
\newcommand{\R}{{\mathbb R}}
\newcommand{\Z}{{\mathbb Z}}
\newcommand{\C}{{\mathbb C}}

\newcommand{\A}{{\mathbb A}}
\newcommand{\cA}{{\mathcal A}}
\newcommand{\cG}{{\mathcal G}}

\newcommand{\LH}{{^L\!H}}
\newcommand{\LG}{{^L\!G}}

\newtheorem{theorem}{Theorem}%[section] (If you want theorem numbered
\newtheorem{lemma}{Lemma}%               with section number.  Same
\newtheorem{proposition}[theorem]{Proposition}
\newtheorem*{principle}{Principle}
\theoremstyle{definition}

\numberwithin{conj}{section}

\newtheorem{definition}{Definition}
\numberwithin{definition}{section}

\numberwithin{question}{section}
\numberwithin{equation}{section}

\theoremstyle{remark}
\newtheorem{remark}{Remark}
\newtheorem*{acknowledgement}{Acknowledgement}
\newtheorem*{Dedication}{Dedication}

\begin{document}

\title[Functorialiy-reduction to the semistable case]{On Langlands 
functoriality- reduction \\ to the semistable case}
\date{}
\author{C.~S.~Rajan}

\address{Tata Institute of Fundamental 
Research, Homi Bhabha Road, Bombay - 400 005, INDIA.}
\email{rajan@math.tifr.res.in}
\dedicatory{To M. S. Raghunathan}

\begin{abstract} We generalize a beautiful method of Blasius and
Ramakrishnan,  that in order to exhibit particular instances
of the Langlands functorial correspondence, it is enough to show that
the correspondence holds in the semistable case, provided the
arithmetic data we are considering is closed with respect to cyclic
base change and descent. 

\end{abstract}

\maketitle

\section{Introduction}
Let $F$ be a global field, and let $H, ~G$ be reductive algebraic
groups defined over $F$. Suppose $G$ is quasi-split, and we have a
morphism 
\[ \phi: \LH \to \LG\]
of the Langlands dual groups $\LH$ and $\LG$ defined over the complex
numbers. The conjectures of Langlands roughly predicts the existence of a
functorial transfer from the collection of automorphic representations
${\mathcal A}(H)$ on $H$ to the collection of automorphic representations
${\mathcal A}(G)$ on $G$.

An example of functoriality is obtained upon specialising $H$ to be
the group consisting of a single element, and $G$ to be $GL_n$ over
$F$.  Let
$G_F$ denote the Galois group of a separable closure $\bar{F}$ of $F$
over $F$, and $\A_F$ denote the adele ring of $F$. The required
functoriality, known as the strong Artin conjecture, is  
the problem of attaching an automorphic representation of $\pi(\rho)$
of $GL_n(\A_F)$
to a linear representation $\rho$ of $G_F$ to $GL_n(\C)$. 
When $n=1$, this correspondence is given  by  the Artin
reciprocity law. The compatiblity
with the local correspondences at all places translates into
\[ L(s, \rho)=L(s, \pi(\rho)),\]
where the $L$-function on the left is the $L$-function defined by
Artin for a Galois representation, and the $L$-function on the right
is the one attached to an automorphic representation of $GL_n(\A_F)$
by the method of Tate-Godement-Jacquet. In particular, the strong
Artin conjecture implies the Artin conjecture  that the functions
$L(s,\rho)$ are entire,  provided $\rho$ does not contain the trivial
representation.  

Another instance of functoriality is obtained by considering a Galois
extension $K/F$. Let $H$ be split over $F$, and let $G=R_{K/F}(H)$,
the Weil restriction of scalars of $H$ considered  over $K$ to $F$. 
At the level of
$L$-groups, there is a diagonal embedding of $\LH$ to $\LG$. The
corresponding lifting is known as base change. When $K/F$ is a cyclic
extension of number fields of prime degree and $H=GL(n)$,   the existence and 
properties of base
change, in particular the descent criterion that invariant
cuspidal automorphic representations lie in the image of base change,
has been shown by
Langlands \cite{L} for $GL(2)$, and by Arthur and Clozel for $GL(n)$
in general \cite{AC}.  

 One of the properties required of the global 
correspondence is that it be  compatible with the local
correspondences at the unramified places as follows: 
given an automorphic representation $\pi$ of
$H(\A_F)$, let $\Pi$ be the corresponding lift to an automorphic
representation of $G(\A_F)$. Let  $v$ be a place of $F$ at which the
local components $\pi_v$  and $\Pi_v$ 
respectively of $\pi$ and $\Pi$  at $v$ are unramified.  The required
compatiblity is that 
\[\phi(S(\pi_v))=S(\Pi_v),\]
where for an unramified representation $\eta$ of a unramified, 
reductive group $H$ over a local field $F$, $S(\eta)$ denotes the
Satake parameter lying in the dual group $\LH$. 

More generally for a local field $F$, let $W_F$ denote the Weil group
of $F$, 
and if $F$ is non-archimedean, let  $W'_F=W_F\times SU_2(\R)$ be the
Weil-Deligne group attached to $F$ (for $F$ archimedean, let
$W'_F=W_F$).  Let $\Phi(G)$ denote the set of
equivalence classes of $L$-parameters, 
\[ \phi: W_F'\to \LG,\]
satisfying the following properties: 
\begin{itemize}
\item   composition with the projection $\LG \to G_F$ is the natural
map occuring in the definition of Weil groups.
\item  $\phi(w)$ is
 semisimple for every $w\in
 W_F$.
\end{itemize}
$\LG$ is the semidirect product of the dual group
$\hat{G}$ defined over the complex numbers and $G_F$, and   the
equivalence relation is taken with respect to the conjugacy action of
$\hat{G}$ on $\LG$. Let $\Pi(G)$ denote the collection of
irreducible, admissible, representations of $G(F)$. 
 Then the local Langlands conjecture predicts a natural bijection
 between $\Phi(G)$, and $L$-packets of representations in
 $\Pi(G)$. The global correspondence is then required to be compatible
 with the local correspondence at all places. 

To enunciate the main principle of the paper, we need a concept of
when the local component of an automorphic representation is  `nice'. 
 \begin{definition}Let $F$ be a  non-archimedean local field, and $G$
be a reductive group defined over $F$. 

\begin{enumerate}
\item $G$ is said to be
{\rm unramified}, if $G$ is quasi-split and splits over an
unramified extension of $F$. 

\item A $L$-parameter $\phi: W'_F\to \LG$  is {\em
semistable}, if $G$ is unramified and if the image of the inertia
group of $F$ with respect to $\phi$ acts
trivially on $\hat{G}$. 
\end{enumerate}
\end{definition}

In other words, a parameter is semistable, if restricted to $W_F$ it
is an unramified parameter in the traditional  sense.  
Since unramified representations are semistable, for almost all places
the local components of an automorphic representation are semistable. 
At the representation theoretic level, a related notion for 
a representation of a split, reductive $p$-adic  group to be
semistable, is 
that the representation
has Iwahori fixed vectors. For $GL(n)$ over local fields, this is
compatible with the definition by parametrizations of the Weil-Deligne
group to the $L$-group $^LGL(n)^0=GL(n)$, given by Bernstein and
Zelevenskii \cite{BZ}. In general, we have the results of
Kazhdan and Lusztig \cite{KL}, 
 relating the different notions of semistability given by 
 the Deligne-Langlands
conjecture.

It was shown by Grothendieck that any $l$-adic representation of
$G_F$ is potentially semistable (see next section). 
Thus the existence of the local Langlands correspondence or a
suitable theory of local base change,  
in partiuclar should imply that given any representation
$\pi\in \Pi(G(F))$ as above, there exists an extension $K$ of $F$,
such that the restriction of the parameter to the Weil group of $K$ is
semistable. 

The main principle we want to illustrate here  is the following:
\begin{principle} \label{principle}
Suppose $F$ is a global field and there is a
morphism 
\[\phi_F: \LH\to \LG\]
of $L$-groups. Assume that the following
conditions are satisfied:
\begin{enumerate}
\item for any extension $K$ of $F$, and a semistable automorphic
representation $\pi$ of $H(\A_K)$, the automorphic lift $\phi(\pi)$ exists. 

\item given $\pi$ an automorphic representation of $H(\A_K)$, 
by  a succession of cyclic base changes of $\pi$, the base changed
automorphic representation is semistable. 

\item the collection of cuspidal automorphic representations $\pi$ of
$G(\A_L)$, satisfy the
descent properties of base change for cyclic extensions
$L/K$ of prime degree,  and $K$ containing $F$, i.e.,
if $\pi$ is any $\text{Gal}(L/K)$-invariant, cuspidal, automorphic
representation of 
$G(\A_L)$, then $\pi$  lies in the image of base change. Further any two
representations of $G(\A_K)$ which base change to $\pi$, differ upto
twisting by an idele class character corresponding to the extension
$L/K$.  
\end{enumerate}

Then the Langlands transfer corresponding to $\phi$,  
can be defined for all cuspidal, automorphic
representations of $H(\A_F)$.

\end{principle}

The principle indicates a more central role for cyclic base change in
the context of establishing reciprocity laws. 
The underlying reasons for  the validity of the  priniciple,  
lies  in  the fact that the local Galois
groups are solvable, and that local extensions can be approximated by
global extensions. The main ingredient that goes into the proof is 
the theory of cyclic base change, and that locally the class of data
considered are potentially semistable.
Using the Grunwald-Wang theorem, it is possible to
produce an infinite family of cylic extensions, over which the extent
of `ramified  behavior' of the given representation comes down, and an
induction argument combined with the properties of base change, allows
us to conclude the theorem. 

The method can be considered as
a refinement of the Chebotarev-Artin trick, and 
appears in   Blasius and Ramakrishnan
\cite{BR}. It has been used in many different contexts, by 
Ramakrishnan  in showing the automorphicity of Rankin-Selberg
$L$-functions  for $GL(2)\times GL(2)$ \cite{R}, and in proving some
cases of Artin's conjecture for Galois representations with solvable
image in $GO(4)$ \cite{R1}. Harris and Taylor have used it to extend
Clozel's results attaching Galois representations to regular, self-dual
automorphic forms on $GL(n)$ \cite{HT}. It is our hope that by placing
the method in a general context, the usefulness and the fundamental
nature of the principle outlined here will become more apparent and known. 

The principal motivation for this paper was  to extend the
geometric correspondence of attaching a $l$-adic representation to
cusp forms on $GL(n)$ over a function field of a curve over a finite
field, provided the cuspidal representation is unramified at all
places. Such a geometric correspondence has been established in the
case $n=2$ by Drinfeld, and by Gaitsgory, Frenkel and Vilonen for
$n\geq 3$ \cite{Lau1}. The results outlined here indicate that it is
indeed possible to obtain the correspondence in general using these
methods, provided the methods of
Gaitsgory, Frenkel and Vilonen can be extended to cover the semistable
case, and if the required properties of base change are established
for $GL(n)$ over function fields. This will then provide a different
and geometric way of obtaining one half of the Langlands
correspondence, than the proof using converse theory and 
 the principle of recurrence given
by Deligne and Piatetskii-Shapiro \cite{Lau}, which reduces the proof
to Lafforgue's theorem attaching automorphic representations to
$l$-adic representations.  

Even in the proof given by  Lafforgue, the
principle of reduction to the semistable case can be applied, and the
proof reduced to the case when the automorphic representation is
semistable, provided cyclic base change for $GL(n)$ over function
fields is available. The advantage is that the compactification of the
space  of shtukas with multiplicity free (reduced) level structure
enjoys nice properties \cite[Theorem III.17]{Laf}, which are not
available in general.

In Section 2, 
we illustrate this principle in the situation of attaching automorphic
representations to $l$-adic representations of $G_F$, $F$ a number
field.  In the last
section we discuss the transfer at the ramified  places, and 
present a method to go from weak to strong lifting
in the context of the methods used in this paper.  Here our attempt is
to clarify the use of different triple product L-functions that appear
in the proof of automorphicity of the Rankin-Selberg convolution for
$GL(2)\times GL(2)$, given by Ramakrishnan \cite{R}.

\begin{acknowledgement} I sincerely thank D. Gaitsgory,
L. Lafforgue and R. Raghunathan for useful discussions and
encouragement. 
\end{acknowledgement}

\begin{Dedication} It is a pleasure to dedicate this paper to
M. S. Raghunathan on his sixtieth birthday. His infectious enthusiasm,
dynamism, attitude towards mathematics and openness have been of great
inspiration to me. We wish him and his family the best in the years to
come. 
\end{Dedication}
 
\section{Reciprocity: reduction to the semistable case}
 Let $F$ be a number field. Denote by
${\cA}(n,F)$ the set of irreducible cuspidal, automorphic
representations of $GL(n, \A_F)$, where $\A_F$ is the adele ring of
$F$. Let  $l$ be  a prime
distinct from $p$, and let $S$ be a finite set of places of $F$
containing the primes of $F$ lying above $l$. Denote by ${\cG}(n,F)$
the set of irreducible $\lambda$-adic  Galois representations 
\[ \rho:G_F\to GL_n(\bar{\Q}_l),\]
 of the absolute Galois group
$G_F$ of an algebraic closure $\bar{F}$ of $F$ over $F$, unramified
outside $S$.  

For a place $v$ of $F$, let $\rho_v: G_v\to 
GL_n(\bar{\Q}_l)$ be the local representation at $v$ associated to
$\rho$, where $G_v$ is the local Galois group. 
\begin{definition} Let $v$ be a finite place of $F$, with residue
characteristic coprime to $l$.  
$\rho_v$ is {\em semistable}  at a place $v$
of $F$, if  $\rho_v|_{I_v}$ is unipotent. $\rho_v$ is {\em potentially
semistable}  at a place $v$
of $F$, if  $\rho_v|_{I_v}$ is quasi-unipotent, i.e., there is a
subgroup $I'_v$ of $I_v$ of finite index, such that $\rho_v|I'_v$ is
unipotent.  

$\rho$ is {\em (potentially) semistable}, if at all finite  places $v$
of $F$, $\rho_v$ is (potentially) 
semistable. (Additionally we can require that $F$ be a totally CM
field). 
\end{definition}

\begin{remark} At the finite places of $F$ lying over $l$, the
corresponding notions of a $l$-adic representation to be semistable
have been defined by Fontaine.
\end{remark}

  It was proved by
Grothendieck that any $\lambda$-adic representation $\rho_v$ is potentially
semistable.  From the
correspondence between $\lambda$-adic representations of the local
Galois group and  representations of the Weil-Deligne group of the
local field as given in \cite[Theorem 4.2.1]{T}, 
we see that the two  definitions of semistability given in the
previous section and the foregoing one are
equivalent. 

Fix an isomorphism of $\bar{\Q}_l$ with $\C$. Let $v$ be a finite
place of $F$ not in $S$. Denote by $\rho(\sigma_v)$ the corresponding 
Frobenius conjugacy class in $GL(n,\C)$. 
\begin{definition} An automorphic representation $\pi$ of $GL_n(\A_F)$
is a (weak) lifting of $\rho$, if at any unramified finite place $v$
of $\rho$ and $\pi$ and $v$ not dividing $l$, we have 
\[ \rho(\sigma_v)=S(\pi_v),\]
where $S(\pi_v)$ is the conjugacy class in  $GL(n,\C)$ defined by the
Satake parameter $S(\pi_v)$. 
\end{definition} 

We illustrate the principle stated above by the following
theorem. If $K$ is a finite extension of $F$, let $\rho^K$ denote the
restriction of $\rho$ to $G_K$. 
\begin{theorem}\label{galois-automorphic}
 Let $\rho: G_F\to GL_n(\bar{\Q}_l)$ be an irreducible
$\lambda$-adic representation of $G_F$. 
Suppose  that  for any finite extension $K$ of $F$, such that the
restriction $\rho^K$ is semistable, there exists a cuspidal, 
automorphic representation $\pi(\rho^K)$ of
$GL(n, \A_K)$ which is a weak lift of $\rho^K$.  
 
Then there exists   a cuspidal,  automorphic representation $\pi(\rho)$
 of\\
 $GL(n, \A_F)$,  which is a weak lifting of $\rho$.  
\end{theorem}

\begin{proof} The proof has three main ingredients. The first
main  ingredient involved in the proof is the structure of  the absolute Galois groups of
local fields.  The results that we need from the structure theory of
local Galois groups,  are 
\begin{itemize}
\item for any place $v$ of a number field $K$, the local Galois group
$G_{K_v}$ is solvable. 

\item{\em Grothendieck's theorem}. Any $\lambda$-adic  representation
is potentially semistable. 
\end{itemize}

Consider the collection of pairs ${\mathcal S}$ of the
form  $(G, p)$, where $G$ is
a solvable group, $p$ is a prime number, 
 and there exists a surjective homomorphism $G\to
\Z/p\Z$. Define the modulus of such a pair to be $(|G|, p)$, and the
prime modulus to be $p$.  
Order these pairs by:
\begin{equation*}
\begin{split}
 (G,p)\leq (G',p')\quad & ~\text{iff}~ |G|\leq |G'|\\
 & \text{or} ~ |G|=|G'| ~\text{and}~ p\leq p'.
\end{split}
\end{equation*} 
Given a $\lambda$-adic representation $\rho: G_F\to GL_n(\bar{\Q}_l)$,
let $S'$ denote the set of finite places $v$ of $F$, at which the
local component $\rho_v$ is not semi-stable. For each place $v\in S'$,
choose a solvable extension $L_v$ of minimal degree over $F_v$ and
with Galois group $H_v$, such that the restriction $\rho_v|_{G_{L_v}}$
is semistable. Consider the collection of pairs  ${\mathcal S}_{\rho}$
given by $(H_v, p)$, satisfying that there exists  a surjective morphism
from $H_v\to
\Z/p\Z$. Define the ramification index $R(\rho)$ of $\rho$ to be the
modulus of a maximal element in ${\mathcal S}_{\rho}$ (depends on the
choice of $L_v$), and let $T$ be the set of places at which the maximum is
attained.

The proof of the theorem is by induction on $R(\rho)$, and we
assume that given a representation $\rho$, 
we have proved the theorem over all extensions  $K$ 
of $F$,  such that $R(\rho^K)<R(\rho)$. \\

The second main ingredient involved in the proof is the Grunwald-Wang
theorem \cite[Chapter XV]{AT}: given a finite set of characters
\[\eta_v :G_v\to \C^*,\]
of order $n_v$  at a finite set of places $v\in T'$ of $F$, there
exists an idele class character $\chi$ of $F$, such that we have an
equality of local components for all $v\in T'$: 
\[ \chi_v=\eta_v.\]
Let $m$ be least common multiple of the orders $n_v$. Let $s$ be the
least integer $\geq 2$, such that $\zeta_{2^s}+\zeta_{2^s}^{-1}\in F$,
but $\zeta_{2^{s+1}}+\zeta_{2^{s+1}}^{-1}\not\in F$, and let 
\[ a_0=(1+\zeta_{2^s})^m.\]
In general, $\chi$ can be chosen to have period $2m$, but the period
can be made equal to $m$, provided  when $F$ is special, the
following auxiliary condition (see \cite[Theorem 5, Chapter
10]{AT}) holds: 
\[ \prod_{v\in T'}\eta_v(a_0)=1. \] 
 In our situation, 
if we assume further that either $\eta_v$ is trivial,  or it is of
order $p$ for  a fixed prime $p$   at all the places $v\in T'$, then
we can arrange for $\chi$ to have order $p$, by adding an extra fixed
finite place $v_0$ to the collection $T'$, and choosing an appropriate
character at $v_0$, such that in the special case appearing in the
hypothesis of the Grunwald-Wang theorem,   the above extra condition
needed for $\chi$ to have order $p$ is satisfied (this is required
only when $F$ is a number field). 

We need a preliminary lemma before applying the Grunwald-Wang theorem:
\begin{lemma} Let $\rho: G_F\to GL_n(\bar{\Q}_l)$ be an irreducible
$\lambda$-adic representation of $G_F$. Then there exists only
finitely many characters $\chi: G_F\to \bar{\Q}_l^*$, such that 
\[\rho\simeq  \rho\otimes \chi.\]
\end{lemma}
\begin{proof}
The ramification of $\chi$ is bounded by the
ramification of $\rho$, and on comparing determinants, we find that
$\chi$ is of order dividing $n$. The lemma follows from the the theorem of
Hermite-Minkowski.  
\end{proof}

We now apply the Grunwald-Wang theorem, and construct characters 
depending on an auxiliary prime $w$ not belonging to $T$. 
Let $T$ be  the set of places of $F$, where the ramification
index of $\rho$ attains it's maximum. For each $v\in T$,  
choose a character $\eta_v$ of order $p$,
where $p$ is the prime modulus and such that the kernel of $\eta_v$ is
a Galois extension of $F_v$ contained in the chosen extension $L_v$. 
Choose,  and fix a place $v_0$ of $F$ not in $T$ (if needed),  
and a character $\eta_{v_0}$, 
such that in the special case the extra condition needed in the
Grunwald-Wang theorem for the global character to have order exactly
the least common multiple of the local orders,  is satisfied (see
\cite[Theorem 5, Chapter 10]{AT}). 
Finally, let $K'$ be the
compositum of the finitely many cyclic extensions of $F$ defined by
the kernels of the characters satisfying the hypothesis of the above
lemma. Choose a prime $w_0\neq v_0$ of $F$ which is inert in $K'$, and let 
\[ \eta_{w_0}=1.\] 
This has the consequence that the fields we construct will be disjoint
from $K'$. Let $w$ be a prime of $F$ distinct from the primes
belonging to $T$ and from $v_0$ and $w_0$, and  
let  $\eta_w=1$. Let $T(w)= T\cup \{w_0\}\cup \{w\}\cup\{v_0\}$. 
By the Grunwald-Wang theorem, we
obtain an idele class  character $\chi^w$,  such
that the local components at $v\in T(w)$ satisfy, 
\[ \chi^w_v=\eta_v.\]
Denote by $K^w$ the corresponding cyclic extension of $F$ defined by
the kernel of $\chi^w$. By construction, $K^w$ satisfies the following
properties:
\begin{itemize}
\item $K^w$ is of order $p$ over $F$, where $p$ is a prime. 

\item $K^w$ is completely split over $F$
at the primes $w_0$ and $w$. In particular, $K^w$  is disjoint from
$K'$. 

\item the ramification index $R(\rho^{K^w})$ is less than
$R(\rho)$. \\

\end{itemize}

The third main ingredient involved in the proof is the theory of base
change, due to Langlands for $n=2$, and by Arthur and Clozel for
general $n$ (\cite{L}, \cite{AC}). Let $K/F$ be a cyclic extension of
prime degree with Galois group
$G(K/F)$. An 
automorphic representation $\Pi$ of  $GL(n, \A_K)$ is invariant with
respect to $G(K/F)$,  if 
$\Pi^{\sigma}\simeq \Pi$, for any $\sigma\in G(K/F)$. Let $\Pi$ be a
cuspidal, automorphic invariant representation of  $GL(n,
\A_{K})$. Then there exists a cuspidal, automorphic representation
$\pi$ of  $GL(n,
\A_{F})$, which base changes to $\Pi$. Further if $\pi_1$ and $\pi_2$
are two such representations which base change to $\Pi$, then 
\[ \pi_1\simeq \pi_2\otimes \psi,\]
where $\psi$ is an idele class character of $F$, corresponding via
class field theory to a character of $G_F\to \C^*$, whose kernel
contains $G_K$. As a consequence of the descent criterion, we have
\begin{lemma} 
Suppose $K_1, \cdots, K_m$ are cyclic, linearly disjoint extensions of
$F$ of prime degree.  For $i=1, \cdots, m$, let 
 $\pi_i$ be a cuspidal,  automorphic representations
of  $GL(n, \A_{K_i})$, invariant
with respect to the actions of $G(K_i/F)$. Assume
further that for any pair of integers $1\leq i, ~j\leq m$, 
the base change $\pi^{ij}$ of $\pi_i$ and $\pi_j$ to the compositum
$GL(n, \A_{K_iK_j})$ are isomorphic. Then there exists a unique, cuspidal
automorphic representation $\pi$ of  $GL(n, \A_{F})$, which base
changes to $\pi_i$ for  $1\leq i \leq m$. 
\end{lemma} 
\begin{proof} We first prove the lemma when $m=2$.
By hypothesis, let $\pi_1', \cdots, \pi_l'$ be the cuspidal
automorphic representations of  $GL(n, \A_{F})$, which base
changes to $\pi_1$, where $l=[K_1:F]$.  
Now the various possible descents differ from one
another by an idele class character on $F$, corresponding to the
extension $K_1/F$. Since $K_2$ is linearly disjoint with $K_1$ over
$F$, these characters remain distinct upon base changing to $K_2$, and
hence the base change $\Pi_i'$ to $K_2$ of the representations $\pi_i'$
remain distinct. 
By hypothesis, the base change of $\Pi_i'$ and $\pi_2$ to the compositum
$K_1K_2$ are isomorphic. Since there are exactly $[K_1K_2:K_2]=l$ such
representations of  $GL(n, \A_{K_2})$, there is precisely one index
$j$ such that $\Pi_j'\simeq \pi_2$, and that gives us the lemma when
$m=2$. 

Now we consider the general case. The descent $\pi_{12}$
living on $F$,  upon base changing to a field $K_l$ different from $K_1$
and $K_2$, satifies the property that upon further base changing to
the fields $K_1K_l$ and $K_2K_l$  becomes isomorphic to $\pi^{1l}$ and
$\pi^{2l}$ respectively. Hence by the uniqueness proved above for two
components, the
base change of $\pi_{12}$ to $K_l$ is isomorphic to $\pi_l$, and we
have the lemma. 
\end{proof}

Now we proceed to the proof of the Theorem.  The restriction\\
$\rho^w:=\rho^{K^w}$ is a $\lambda$-adic representation of $G_{K^w}$,
with ramification index less than that of $\rho$. By induction hypothesis,
we have a cuspidal,  automorphic representations $\pi^w$
of  $GL(n, \A_{K^w})$, which is a (weak) lifting of $\rho^w$. By
strong multiplicity one, $\pi^w$ is invariant with respect to
the action of $G(K^w/F)$,   since $\rho^w$ is
invariant. 

Hence by the lemma proved above, we have a unique cuspidal,
automorphic representation $\pi$ of $GL(n, \A_{F})$, which base
changes to $\pi^w$ for each $w$. By the inductive hypothesis, 
the local components  of $\pi^w$ and $\rho^w$ at an unramified  prime of
$K^w$ lying above $w$ agree. But $w$ splits completely in $K^w$, and
thus these local components are isomorphic respectively to the local
components of $\rho$ and $\pi$ at the place $w$. Hence 
at any place $w$ of $F$,  not belonging to $T\cup\{v_0\}$ and where
$\rho_w$ and $\pi_w$ are unramified,  we obtain that the Frobenius
conjugacy class defined by $\rho_w$ is equal to the Satake parameter
associated to $\pi_w$. By choosing a different $v_0$, and appealing to
strong multiplicity one, we obtain the compatiblity of the lift at
$v_0$ too, and that proves the theorem.
\end{proof}

\begin{remark} The above theorem can be used in the reverse direction
too, from knowing how to attach Galois representations to `semistable,
cuspidal,  automorphic representations' to attaching Galois representation to 
any cuspidal, automorphic representation. For this, one has to know
that locally the representations are potentially semistable, either
via the local Langlands correspondence or using the notion that a
 representation being semistable amounts to the representation having
non-zero  
 Iwahori fixed vectors. For $GL(n)$ over local fields, this is
compatible with the definition by parametrizations of the Weil-Deligne
group to the $L$-group \cite{BZ}, and  in general we have the results of
Kazhdan and Lusztig \cite{KL} establishing the Deligne-Langlands
conjecture.  
 
\end{remark}

\subsection{Function fields} 
Let $F$ be a global field of positive characteristic. 
The Langlands correspondence for $GL(2)$ was established by Drinfeld. 
In \cite{Laf}, Lafforgue has established the Langlands conjecture, giving a canonical
bijection between the collection of cuspidal, automorphic
representations of $GL_n(\A_F)$,  and the collection of irreducible
$\lambda$-adic representations  of $G_F \to GL_n(\bar{\Q}_l)$
unramified outside a finite set of
places of $F$,  compatible with the
local Langlands correspondence at all places $v$ of $F$. 
The method of Lafforgue is to attach  a Galois representation
$\rho(\pi)$  to
any cuspidal automorphic representation  $\pi\in {\cA}(n,F)$. The reverse
map $\rho\mapsto \pi(\rho)$,  is the `principle of recurrence',  due to
Deligne, Laumon and Piatetskii-Shapiro. This 
 is given by appealing to the converse theorems proved by
Cogdell and Piatetskii-Shapiro, and the results of Grothendieck 
that the $L$-functions
attached to Galois representations have nice analytic properties. 

Using geometric methods, Gaitsgory, Frenkel and Vilonen have
 established the correspondence 
$\rho\mapsto \pi(\rho)$, attaching 
a cuspidal automorphic representation $\pi(\rho)$ to any unramified
Galois representation $\rho$ of $G_F$, and compatible with the local
correspondence (see \cite{Lau1}).   The main motivation for this paper
is to show that this correspondence can be extended to cover all
representations $\rho\in {\cG}(n,F)$, provided we know the
correspondence for all `semistable representations', and we have the
Langlands-Arthur-Clozel theory of cyclic base change for extensions of
prime degree. 

Although the results that we require about base change
follow from Lafforgue's results, it would be desirable to have a
direct proof of base change, and to characterize the local lifts in
terms of character identities, since it reduces the proof of
Lafforgue's theorem to the case when the automorphic representation is
semistable.  Granting  the existence of local base
change asserting that the local components are potentially
semistable, by the   application of the principle given above, 
we only  need to attach Galois representations to
semistable automorphic representations. This  amounts to
considering the space of shtukas with reduced  level
structure and their  compactifications, which 
enjoys nice properties \cite[Theorem III.17]{Laf}, and  are not
available in general.

\section{Weak to strong lifting: Matching of $L$ and $\epsilon$-factors}
Our main aim in this section  is to indicate how to go from a weak
lifting to a strong lifting, where we obtain  the matching of  $L$ and
$\epsilon$-factors match at all places, assuming that in the case of
lifting of semistable data, the local $L$ and $\epsilon$-factors match
at all places. The applications we have in mind here are  to clarify
the use of different $L$-functions associated to triple products that
appears in the proof of automorphicity of the Rankin-Selberg
convolution for $GL(2)\times GL(2)$ established by Ramakrishnan
\cite{R}, and to provide a slightly different perspective of  Laumon's theorem
proving the factorisation of the global epsilon factor occuring in 
Grothendieck's functional equation for $L(s,\rho)$,  in terms of the
local constants constructed by Deligne and Langlands. 

 Again the main ingredient is the theory of base change,
combined together with Grunwald-Wang theorem, which allows us to
concentrate  our attention at only one place. The additional
assumption that we require, is the existence of  Jacquet-Shalika type 
bounds at the semistable places and the global functional equation. 
For the case of automorphic representations on $GL(n)$,
Jacquet-Shalika type bounds are available from the existence and
properties of Rankin-Selberg convolution proved by Jacquet,
Piatetskii-Shapiro, Shalika and Shahidi
 (\cite{JS}, \cite{JPSh}, \cite{Sh}).  

For the mathching of epsilon factors, we
assume the factorisation of the global $\epsilon$-factor, and then use
the defining properties of the epsilon factors, together with the 
functional equation. We illustrate the matching at all places of the
local $L$ and $\epsilon$-factors, in the context of the principle.  

We first abstract a notion of an `arithmetic $L$-data' $\cA$,
which is stable with respect to base change.

\begin{definition}
An {\em arithmetic $L$-data}, is a collection $\cA$ satisfying  the following:
\begin{enumerate}
\item to any element $\pi\in \cA$, we have an associated global field
$K$. We will say that $\pi$ is defined over $K$.  
Further there is an $L$-function, 
\[ L(s, \pi)=\prod_{v\in \Sigma_K}L_v(s,\pi),\]
indexed over the places $v\in \Sigma_K$ of $K$. It is required that
$L_v(s,\pi)$ are meromorphic and non-vanishing on the entire plane. 
  At the archimedean
places, $L_v(s,\pi)$ can be expressed in terms of
$\Gamma$-functions. At the finite places, we require that 
\[ L_v(s,\pi)=\prod_{i=1}^{d_v}(1-\alpha_{v,i}q_v^{-s}),\]
where $q_v$ is the cardinality of the residue field, and this will be
referred to as the local $L$-factors of $\pi$.  Further it is
required that the product defining $L(s,\pi)$ be absolutely convergent
in some right half plane.

\item{\em Functional equation.}
$L(s,\pi)$ has a meromorphic analytic
continuation to the entire plane, and satisfies a functional equation
in the form, 
\[ L(1-s,\tilde{\pi})=\epsilon(s,\pi)L(s,\pi),\]
where $\tilde{\pi}\in \cA$  is the
`contragredient' of $\pi$, and is defined over $K$. We have that
$\tilde{\tilde{\pi}}=\pi$, i.e., the contragredient operation is an
involution on $\cA$. $\epsilon(s,\pi)$,
called the global $\epsilon$-factor associated to $\pi$, is an 
entire, non-vanishing function of $s$.  
Define the archimedean component $L_{\infty}(s,\pi)$ of $L(s,\pi)$ by
\[ L_{\infty}(s,\pi)=\prod_{v\in \Sigma_{\infty}}L_v(s,\pi),\]
the product taken over the collection of archimedean places
$\Sigma_{\infty}$ of $K$. 

\item{\em Convolution with Dirichlet characters.} If $\pi$ is defined
over $K$, and $\chi$ is a Dirichlet character of $K$, i.e., an idele
class character of finite order, then there is an element $\pi\otimes
\chi\in \cA$. With the above notation for $L_v(s,\pi)$, 
at a finite place $v$ of $K$ where $\chi$ is
unramified, the local $L$-factor is given by  
\[ L_v(s, \pi\otimes
\chi)=\prod_{i=1}^{d_v(\pi)}(1-\chi(\omega_v)\alpha_{v,i}(\pi)q_v^{-s})^{-1},\]
where $\omega_v$ is an uniformising parameter for $K_v$.
The contragredient of $\pi\otimes \chi$ is $\tilde{\pi}\otimes
\chi^{-1}$. 

\item{\em Base change.}
For any cyclic extension $L/K$ with Galois group $G(L/K)$ and $\pi$ defined over $K$, 
we have a base change element $B_{L/K}(\pi)\in \cA$ defined over $L$,
with the local components of the $L$-function at any place $w$ of $L$
dividing a place $v$ of $K$ given by,
\[ L_w(s,B_{L/K}(\pi))=\prod_{\chi\in 
\widehat{G(L/K)}}L_v(s,\pi\otimes \chi),\]
where $\widehat{G(L/K)}$ denotes the set of characters of $G(L/K)$. Here
we are using abelian class field theory in thinking of $\chi$ also as  an
idele class character, with abuse of notation. If $v$ is a place of
$K$, and $w$ a place of $L$ of degree one over $K$, then
\[L_v(s,\pi)=L_v(s,B_{L/K}(\pi)).\] 

\item{\em Potential semistability.} There is a concept of $\pi$ being
semistable at any place of $K$, satisfying:
\begin{itemize}
\item $\pi$ is semistable at all but
finitely many places of $K$. 

\item if $\pi$ is semistable at $v$, and $L/K$ is a cyclic extension,
then $B_{L/K}(\pi)$ is semistable at all places of $L$ lying over
$v$. 

\item given a place $v$ of $K$, 
there is sequence of cyclic extensions $L=K_n\supset
K_{n-1}\supset \cdots \supset K_0=K$, with $K_{i+1}/K_i$ cyclic such
that 
$v$ splits completely in $L$, and 
the successive base changes of $\pi$ defined with respect to
this filtration exist, with the property that 
\[ B_{L/K}(\pi):=B_{L/K_{n-1}}\circ\cdots \circ B_{K_1/K}(\pi),\]
is semistable at all places of $L$ not dividing $v$. 
\end{itemize}

\end{enumerate}
\end{definition}
\begin{remark} It is possible to relax the condition of $\pi$ being potentially
semistable at all places of $K$, and instead require that the above
conditions be valid for a particular class of places of the global
fields (for example the class of totally real fields), 
closed with respect to divisiblity. 

The last condition in the notion of potential semistability, is
motivated by the expectation that   locally the components of
automorphic representatios are potentially semistable, and the use of
Grunwald-Wang theorem as in the proof of Theorem \ref{galois-automorphic}. 
We have stated it in this form, rather than put in a condition of
local base change and declaring that the local components are
potentially semistable. 
\end{remark} 

\begin{definition} Let $\cA_1, ~\cA_2$ be  two classes of arithmetical
$L$-data as above. Given $\pi_1\in \cA_1$ and $\pi_2\in \cA_2$, we say
that $\pi_1$ and $\pi_2$ are {\em weakly compatible} (or weakly match, or is
a weak lift of one..), if at all places $v$ of $K$, where $\pi_1$ and
$\pi_2$ are semistable, we have
\[ L_v(s,\pi_1)=L_v(s,\pi_2).\] 

$\pi_1$ and $\pi_2$ are {\em strongly compatible} if at all  places $v$
of $K$, we have $L_v(s,\pi_1)=L_v(s,\pi_2)$. 

If $\pi_1$ and $\pi_2$ are
strongly compatible, we have as a consequence that
$\epsilon(s,\pi_1)=\epsilon(s,\pi_2)$. 

Let $R: \cA_1\to \cA_2$ be a map. We say that $R$ is a {\em weak
(resp. strong) transfer} (or weak (resp. strong) lift, or weakly
(resp. strongly)
compatible), if for every $\pi_1\in \cA_1$, $\pi_1$ and $R(\pi_1)$ are
weakly (resp. strongly) compatible.
\end{definition}

\begin{definition} $\pi\in \cA$ satisfies the {\em Jacquet-Shalika
(JS)} condition at a place $v$ of $K$, if $L_v(s,\pi)$ is holomorphic
in the half plane $\text{Re}(s)\geq 1/2$. 
\end{definition}

\begin{remark} As remarked above, this condition is satisfied for
isobaric, automorphic representations of $GL(n, \A_F)$, by the results
of Jacquet and Shalika \cite{JS}. 
\end{remark}

\begin{proposition} \label{weaktostrong}
Suppose  we are given two sets of arithmetic
data $\cA_1$ and $\cA_2$ as above,  and a weak transfer $R:\cA_1\to\cA_2$. Assume
that either $\cA_1$ or $\cA_2$ has the property, that for any element
$\pi$ in one of them, and any place $v$ of $K$ at which $\pi$ is
semistable, $L(s,\pi)$ satisfies the Jacquet-Shalika
property at $v$.  

Then $R$ is a strong transfer. 

\end{proposition}
\begin{proof} Given $\pi\in \cA_1$ and a finite  place $v$ of $K$, we can find
a  sequence of cyclic extensions $L=K_n\supset
K_{n-1}\supset \cdots \supset K_0=K$, with $K_{i+1}/K_i$ cyclic,
 such that $v$ splits completely in $L$, and 
\[B_{L/K}(\pi), ~B_{L/K}(\tilde{\pi}), ~B_{L/K}(R(\pi)),
~~B_{L/K}(\widetilde{R(\pi)})\] 
are  semistable at all places of $L$ not lying over $v$.
Since the local $L$-factors do not change for primes splitting
completely in $L$, 
it is enough now to prove the equality of the $L$-factors at a place
$w$ of $L$  lying above $v$. 

Let $w'$ be a place of $L$ not lying over
$v$. Again we can assume that there is a sequence 
 cyclic extensions $M=L_m\supset
L_{m-1}\supset \cdots \supset L_0=L$, with $L_{j+1}/L_j$ cyclic,
such that 
\[B_{M/K}(\pi), ~B_{M/K}(\tilde{\pi}), ~B_{M/K}(R(\pi)), 
~B_{M/K}(\widetilde{R(\pi)})\]
 are
semistable at all places of $M$. 

We argue by induction on $j$
ranging from $1$ to $m$.  Suppose $M_1/M_2$ is a cyclic extension, and
we have $\pi$ defined over $M_2$. Assume further that
$B_{M_1/M_2}(\pi)$ exista and satisfies the  Jacquet-Shalika condition
at a place $w$ of $M_1$ lying over a place $v$ of $M_2$. 
Then from the equality, 
\[ L_w(s,B_{M_1/M_2}(\pi))=\prod_{\chi\in
\hat{G(M_1/M_2)}}L_v(s,\pi\otimes \chi),\] 
we obtain that $L_v(s,\pi)$ also satisfied the {\em JS}-condition,
since we have the local $L$-factors are non-vanishing everywhere. 
Hence,  we obtain by  induction that the local
components of $B_{L/K}(\pi)$ and $B_{L/K}(R(\pi))$ also satisfy the
Jacquet-Shalika condition at any place of $L$. 

We further have the functonal equations, 
\[ \begin{split}
L(1-s,\tilde{\pi})& =\epsilon(s,\pi)L(s,\pi)\\
\text{and}\quad  L(1-s,~\widetilde{R(\pi})) & =\epsilon(s,\pi)L(s,R(\pi)).
\end{split}
\]
Dividing these two functional equations, and using the fact that the
local components are same at all the places not lying over $v$ by the
hypothesis and  semistability condition, we obtain 
\[
\frac{L_w(1-s,\tilde{\pi})}{L(1-s,\widetilde{R(\pi)})}=\frac{\epsilon(s,{\pi})}
{\epsilon(s,R(\pi))}\frac{L_w(s,{\pi})}{L_w(s,{R(\pi)})}.\]
From the Jacquet-Shalika condition  we have that the right hand side
is holomorphic and  non-vanishing in the half plane $\text{Re}(s)\geq
1/2$. From the above
functional equation and the assumptions about
$\epsilon$-factors,, using the fact that contragedient is an
involution, we obtain that 
\[\frac{L_w(s,{\pi})}{L_w(s,{R(\pi)})}=
\frac{\prod_{i=1}^{d_v(R(\pi))}(1-\alpha_{v,i}(R(\pi))q_v^{-s})}{\prod_{i=1}^{d_v(\pi)}(1-\alpha_{v,i}(\pi)q_v^{-s})},\] is an 
entire, non-vanishing function of $s$.  But then the right hand side
has to be a constant. 

The proof of mathcing of $L$-factors at the archimedean places, is
similar, where we observe that a ratio of products of Gamma
functions cannot be entire unless it is a constant.  
\end{proof}

\begin{remark} In the context of the applicablity of converse
theorems, a similar statement is proved in \cite[Proposition
4.1]{BR}. Instead of
cyclic base change, what is required in \cite{BR},  are converse
theorems and strong multiplicity one.  Indeed by the converse theorems
we have attached two possible (cuspidal) automorphic representation,
which agree at all places except the given ones, and by strong
multiplicity one, we have the matching at all places. 
\end{remark}

\subsection{Automorphicity of Rankin-Selberg.} 
The principle outlined in this paper, has been derived from  the
 method of proof of
 the automorphicity of the Rankin-Selberg convolution for
$GL(2)\times GL(2)$ proved by Ramakrishnan in \cite{R}. Here we
clarify the use of different triple product $L$-functions appearing in
 Ramakrishnan's proof, and show that it is enough to consider the triple product $L$-function
 constructed by Shahidi. 

By the principle outlined in the first part, to prove the automorphicity of the Rankin-Selberg convolution, it is enough
to prove it  in the case when the automorphic
representations involved have semistable reduction at all places. To
prove this for $GL(2)\times GL(2)$, 
we have to show the following:  if $\pi_1$ and $\pi_2$
are semistable,  cuspidal,  automorphic representations of $GL(2,
\A_K)$ for a number field $K$, and if $\eta$ is an arbitrary
cuspidal, automorphic representation of $GL(2,\A_K)$ with the set of
places of $K$ at which $\eta$ is ramified disjoint from the set of
places at which either $\pi_1$ or $\pi_2$ is ramified, then the triple
product $L$-function $L(s,\pi_1\times \pi_2\times \eta)$ is `nice' in
the sense of converse theory, i.e., it has an analytic continuation to
the entire plane, is entire and 
satisfies a suitable functional equation, and is
bounded in vertical strips. 

One knows by the work of Shahidi
\cite{Sh2}, and Kim and Shahidi \cite{KS}, that the triple product
$L$-function can be analytically continued as an entire function, and
that it satisfies a suitable functional equation. Further by
Ramakrishnan \cite{R},  and more generally by Gelbart-Shahidi
\cite{GS}, these $L$-functions are bounded in vertical
strips. The required matching at the semistable places of the $L$ and
$\epsilon$-factors is given by Part i) of Theorem 3.5 in \cite{Sh3} at
the non-archimedean places, and for the archimedean places it is given
by the archimedean Langlands conjecture proved by Langlands. 

This gives us the weak lifting for semistable, cuspidal, automorphic
representations. It can be seen along the lines of the proof
of Theorem \ref{galois-automorphic},  that the principle is valid in
this setting, as the necessary hypothesis are satisfied, and hence we
obtain a weak lifting in general. By   Proposition \ref{weaktostrong}, we
obtain a strong lifting too. 

\begin{remark}
In view of the  results on lifting from classical groups to $GL_n$, it
seems possible to establish the lifting for semistable representations
(those with Iwahori fixed vectors) from classical groups to $GL_n$, 
by extending Proposition 3.1 of \cite{CKPS} to  the general case. 
\end{remark}

\subsection{Function fields and Laumon's factorisation theorem.} 
We now consider the situation over function fields. For
the rest of the discussion, we are assuming the existence and
properties of cyclic base change for $GL(n)$ as outlined
above. In order to attach an automorphic representation to a Galois
representation,  
it is enough now to apply the `principle of recurrence' of Deligne and
Piatetski-Shapiro and obtain a weak automorphic lifting of a semistable,
irreducible   Galois representation $\rho$ of degree $n$.  
The analytic properties of the $L$-function attached to a
$\lambda$-adic representation has been proved by Grothendieck. 
To apply the converse theorem, it remains  to
provide a factorisation of the epsilon factor 
$\epsilon(s,\rho\otimes \rho(\eta))$ occuring in Grothendieck's
functional equation, in terms of the local constans constructed by
Deligne and Langlands, as given by the theorem of  
Laumon \cite[Section 3]{Lau}. Here    $\rho(\eta)$ is the Galois
representation of degree $m$ (less than $n$)  
attached to a cuspidal, automorphic representation of $GL(m, \A_K)$ by
Lafforgue \cite{Laf}. The factorisation of the $\epsilon$-constant for
$\rho(\eta)$ follows from Lafforgue's theorem, or from Lafforgue's
theorem in the semistable case, and applying Theorem
\ref{galois-automorphic} and Proposition \ref{weaktostrong}. 
We can further assume that the ramified
primes of $\eta$ are disjoint from the ramified primes of
$\rho$. Thus we require  Laumon's theorem in the special case
when $\rho$ is semistable, and assuming the factorisation for $\rho$ 
and $\rho(\eta)$ (with disjoint ramificiation loci),  to prove the
factorisation for the    tensor product $\rho\otimes \rho(\eta)$. 
The converse theorem machinery proves the
existence of a weak automorphic lift for $\rho$.
As a consequence,  we also obtain Laumon's theorem
proving the factorisation of the global epsilon factor occuring in 
Grothendieck's functional equation for $L(s,\rho)$,  in terms of the
local constants constructed by Deligne and Langlands.

\end{document}